\documentclass[9pt,conference]{IEEEtran}

\usepackage{fancyhdr}
\usepackage{caption}
\usepackage{subcaption}
\usepackage{framed}
\usepackage[normalem]{ulem}
\usepackage{cite}
\usepackage{amsmath}
\usepackage{amsthm}
\usepackage{amssymb}
\usepackage{amsfonts}
\usepackage{enumerate}
\usepackage{comment}
\usepackage{csquotes}

\newcommand{\complex}{\mathbb{C}^L}

\ifCLASSINFOpdf
   \usepackage[pdftex]{graphicx}
\else
   \usepackage[dvips]{graphicx}
\fi
\usepackage[colorlinks=true,linkcolor=blue]{hyperref}
\newtheorem{theorem}{Theorem}[section]

\newtheorem{definition}[theorem]{Definition}

\makeatletter
 \let\old@ps@headings\ps@headings
 \let\old@ps@IEEEtitlepagestyle\ps@IEEEtitlepagestyle
 \def\confheader#1{%
 \def\ps@headings{%
 \old@ps@headings%
 \def\@oddhead{\strut\hfill#1\hfill\strut}%
 \def\@evenhead{\strut\hfill#1\hfill\strut}%
 }%
 \def\ps@IEEEtitlepagestyle{%
 \old@ps@IEEEtitlepagestyle%
 \def\@oddhead{\strut\hfill#1\hfill\strut}%
 \def\@evenhead{\strut\hfill#1\hfill\strut}%
 }%
 \ps@headings%
 }
 \makeatother

\confheader{%
2021 Online International Conference on
Computational Harmonic Analysis (Online-ICCHA2021)
 }

\begin{document}
%\fancyhdr[C]{2016 IEEE 24th International Requirements Engineering Conference} %zentrierte Kopfzeile

%\makeatletter
%\newcommand*\titleheader[1]{\gdef\@titleheader{#1}}
%\AtBeginDocument{%
%  \let\st@red@title\@title
%  \def\@title{%
%    \bgroup\normalfont\large\centering\@titleheader\par\egroup
%    \vskip1.5em\st@red@title}
%}
%\makeatother

%
% paper title
% can use linebreaks \\ within to get better formatting as desired
\title{Spark Deficient Gabor Frames for Inverse Problems}
%\titleheader{2016 IEEE 24th International Requirements Engineering Conference}

% author names and affiliations
% use a multiple column layout for up to three different
% affiliations
% \author{\IEEEauthorblockN{Vasiliki Kouni\\
% Chair for Mathematics of Information Processing,
% RWTH Aachen University, Aachen, Germany\\
% kouni@mathc.rwth-aachen.de}
% \and \IEEEauthorblockN{Holger Rauhut\\
% Department of Somewhere\\
% email@university.com}
% \and \IEEEauthorblockN{Author 3\\
% Department of Somewhere\\
% email@university.com}
% \and \IEEEauthorblockN{Author 4\\
% Department of Somewhere\\
% email@university.com}}

% conference papers do not typically use \thanks and this command
% is locked out in conference mode. If really needed, such as for
% the acknowledgment of grants, issue a \IEEEoverridecommandlockouts
% after \documentclass

% for over three affiliations, or if they all won't fit within the width
% of the page, use this alternative format:
% 
\author{\IEEEauthorblockN{Vasiliki Kouni\IEEEauthorrefmark{1}\IEEEauthorrefmark{2} and
Holger Rauhut\IEEEauthorrefmark{1}}
\IEEEauthorblockA{
\IEEEauthorrefmark{1}
RWTH Aachen University, Aachen, Germany\\
%Atlanta, Georgia 30332--0250\\ 
Email: kouni@mathc.rwth-aachen.de, rauhut@mathc.rwth-aachen.de}
\IEEEauthorblockA{
\IEEEauthorrefmark{2}
National and Kapodistrian University of Athens, Athens, Greece\\
Email: vicky-kouni@di.uoa.gr}}

% use for special paper notices
%\IEEEspecialpapernotice{(Invited Paper)}

% make the title area
\maketitle

% For peer review papers, you can put extra information on the cover
% page as needed:
% \ifCLASSOPTIONpeerreview
% \begin{center} \bfseries EDICS Category: 3-BBND \end{center}
% \fi
%
% For peerreview papers, this IEEEtran command inserts a page break and
% creates the second title. It will be ignored for other modes.
\IEEEpeerreviewmaketitle

\begin{abstract}
    In this paper, we apply star-Digital Gabor Transform in analysis Compressed Sensing and speech denoising. Based on assumptions on the ambient dimension, we produce a window vector that generates a spark deficient Gabor frame with many linear dependencies among its elements. We conduct computational experiments on both synthetic and real-world signals, using as baseline three Gabor transforms generated by state-of-the-art window vectors and compare their performance to star-Gabor transform. Results show that the proposed star-Gabor transform outperforms all others in all signal cases.
\end{abstract}

\section{Introduction}
% no \IEEEPARstart
% Linear inverse problems deal with recovering a signal $x\in\mathbb{V}^L$ ($\mathbb{V}=\mathbb{R}$ or $\mathbb{C}$) from a set of corrupted observations $y=Ax+e\in\mathbb{V}^K$ ($K\leq L$) and they usually are ill-posed.
We address two ill-posed inverse problems: Compressed Sensing (CS) and speech denoising. To remedy their ill-posedness, we assume $x$ is analysis sparse \cite{rauhut}. The optimization problems occuring under this scenario are --for CS and denoising respectively-- the following:
\begin{align}
    \label{regl1}
    &\min_{x\in\mathbb{V}^L}\|\Phi x\|_1
    %+\frac{\mu}{2}\|x-x_0\|_2^2
    \quad\text{subject to}\quad \|Ax-y\|_2\leq\eta\\
    \label{regl1d}
    &\min_{x\in\mathbb{V}^L}\|\Phi x\|_1\quad\text{subject to}\quad \|x-y\|_2\leq\eta,
\end{align}
$A\in\mathbb{V}^{K\times L}$ ($\mathbb{V}=\mathbb{R}$ or $\mathbb{C}$, $K<L$).
We turn to analysis sparsity instead of its synthesis twin \cite{figu} due to some advantages the former has, e.g. lower computational cost, the optimization algorithm used may need less measurements for perfect reconstruction etc.

\subsection{Motivation and key contributions}
We are motivated by related works proposing either analysis operators $\Phi$ associated to  full spark frames or a finite difference operator \cite{elad,genzel,rauhut}. In a similar spirit, we employ spark deficient Gabor frames (SDGF). Our main contributions are the following: a) we generate an SDGF, associate to it a Gabor analysis operator/digital Gabor transform (DGT) called star-DGT and use it as a sparsifying transform in both analysis CS and denoising b) we compare numerically star-DGT with three other DGTs,
based on famous windows of time-frequency analysis, on synthetic and real-world data. Our experiments show that our method outperforms all others, consistently for all signals, in both CS and denoising.

\section{main results}
\begin{definition}
A \textit{discrete Gabor frame} $(g,a,b)$ \cite{dgs} is a collection of time-frequency shifts of a window vector $g\in\mathbb{C}^L$, expressed as
\begin{align}
    \label{gaborsystem}
    g_{n,m}(l)&=e^{2\pi imbl/L}g(l-na),\quad l\in[L],
\end{align}
which spans $\mathbb{C}^L$ \cite{sparkmal}. Here, $a,\,b$ denote time and frequency parameters respectively, $n\in[N]$ with $N=L/a\in\mathbb{N}$ and $m\in[M]$ with $M=L/b\in\mathbb{N}$ denote time and frequency shift indices respectively.
\end{definition}
The number of elements in $(g,a,b)$ according to \eqref{gaborsystem} is $P=MN=L^2/ab$ and since $(g,a,b)$ is a frame, we have $ab<L$.
\begin{definition}
The \textit{Gabor analysis operator} associated to a Gabor frame is defined as
\begin{equation}\label{gabcoeff}
    \Phi_g:\complex\mapsto\mathbb{C}^{M\times N}:x\mapsto\sum_{l=0}^{L-1}x_l\overline{g(l-na)}e^{-2\pi imbl/L},
\end{equation} for $m\in[M],\,n\in[N]$.
\end{definition}
\begin{definition}[\cite{dang}]
To the \textit{Zauner} matrix
$\mathcal{Z}=\big(\begin{smallmatrix}
0 & \beta\\
1 & \beta
\end{smallmatrix}\big)$, $\beta=L-1$, $L\in\mathbb{Z}$,
corresponds the unitary $U_\mathcal{Z}$ given by the explicit formula \cite{dang}
\begin{equation}\label{sic}
    U_\mathcal{Z}=\frac{e^{i\theta}}{\sqrt{L}}\sum_{u,v=1}^L\tau^{\beta^{-1}(\beta(u-1)^2-2(u-1)(v-1))}e_ue_v,
\end{equation}
with $\theta$ an arbitrary phase, $\beta\beta^{-1}\equiv1\mathrm{modulo}L$, $\tau=-e^{\frac{i\pi}{L}}$ and $e_u,\,e_v$ the standard basis vectors.
\end{definition}
\begin{definition}[\cite{sparkmal}]
The \textit{spark} of a set $F$ of $P$ vectors in $\complex$ is the size of the smallest linearly dependent subset of $F$. A frame $F$ is full spark if and only if every set of $L$ of its elements is a basis, otherwise it is spark deficient.
\end{definition}
\begin{theorem}[\cite{dang}] \label{sdgf}
Let $L\in\mathbb{Z}$ such that $2\nmid L$, $3\mathrel{|} L$ and $L$ is square-free. Then, any eigenvector of the Zauner unitary matrix $U_\mathcal{Z}$ generates a spark deficient Gabor frame for $\complex$.
\end{theorem}
In order to produce an SDGF and apply its associated analysis operator in \eqref{regl1} and \eqref{regl1d}, we first choose an ambient dimension $L$ that fits the assumptions of Theorem \ref{sdgf}. Then, we perform the spectral decomposition of $U_\mathcal{Z}$ in order to acquire its eigenvectors. Since all the eigenvectors of $U_\mathcal{Z}$ generate SDGFs, we may choose an arbitrary one, call it \textit{star window} from now on and denote it as $g_*$. We call the analysis operator associated with such an SDGF \textit{star-DGT} and denote it $\Phi_{g_*}$.
\section{Numerical Experiments}
We solve \eqref{regl1} and \eqref{regl1d}, for synthetic (complex and real-valued) and real-world speech signals respectively, using star-DGT along with three other DGTs emerging from famous window vectors (Gaussian, Hann, Hamming or itersine). Fig. \ref{} presented in the next page, show a) for CS, the 4 relative reconstruction error decays as the number of measurements increases b) for speech denoising, how the 4 MSEs scale as the noise's standard deviation increases. In all cases, star-DGT (blue line) outperforms the rest of DGTs, consistently for all signals and for different choices of ambient dimension with time-frequency parameters for each signal.
\section{Conclusion}
In the present paper, we took advantage of a window vector to generate a spark deficient Gabor frame and introduced a (highly) redundant Gabor transform, i.e. the star-DGT, associated with this SDGF. We then applied star-DGT to analysis Compressed Sensing and analysis-sparsity-based speech  denoising, along with three other DGTs generated by state-of-the-art window vectors in the field of Gabor Analysis. Our experiments confirm improved performance: star-DGT outperforms all others for both synthetic and real-world data. Future directions will include the combination of deep learning architectures with star-DGT, as well as the extension of the presented framework to largescale problems.

\begin{figure}
  \subfloat[]{%
  \begin{minipage}{\linewidth}
  \centering
  \includegraphics[width=0.8\textwidth]{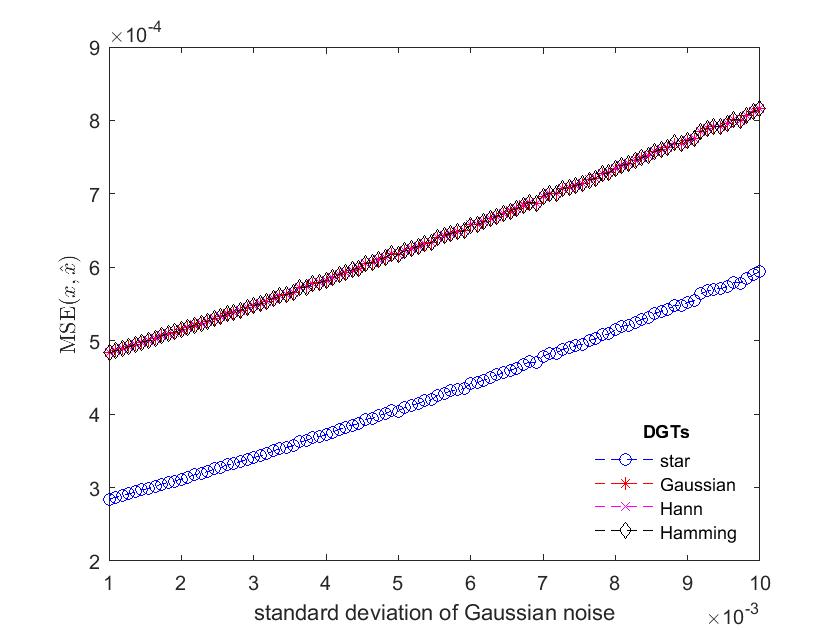}\hfill
  \includegraphics[width=0.8\textwidth]{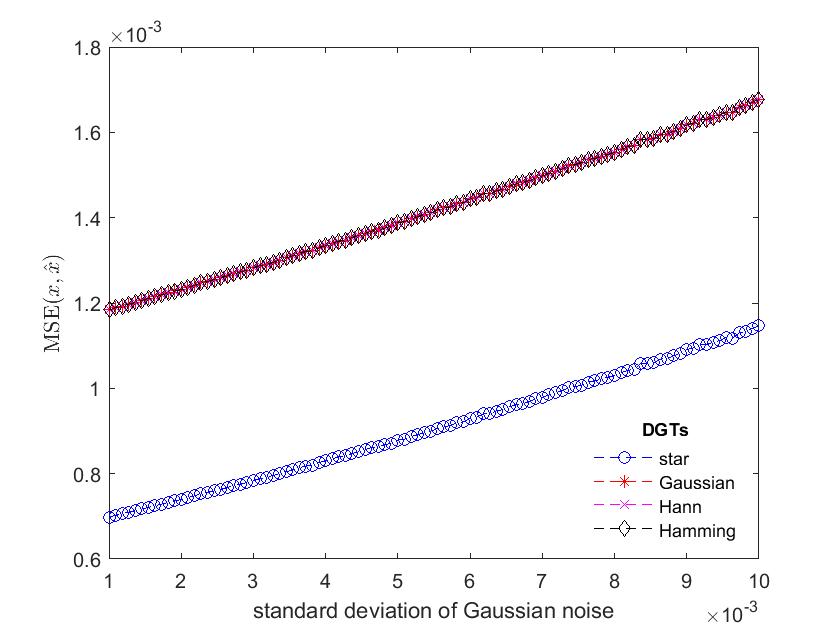}\hfill
  \includegraphics[width=0.8\textwidth]{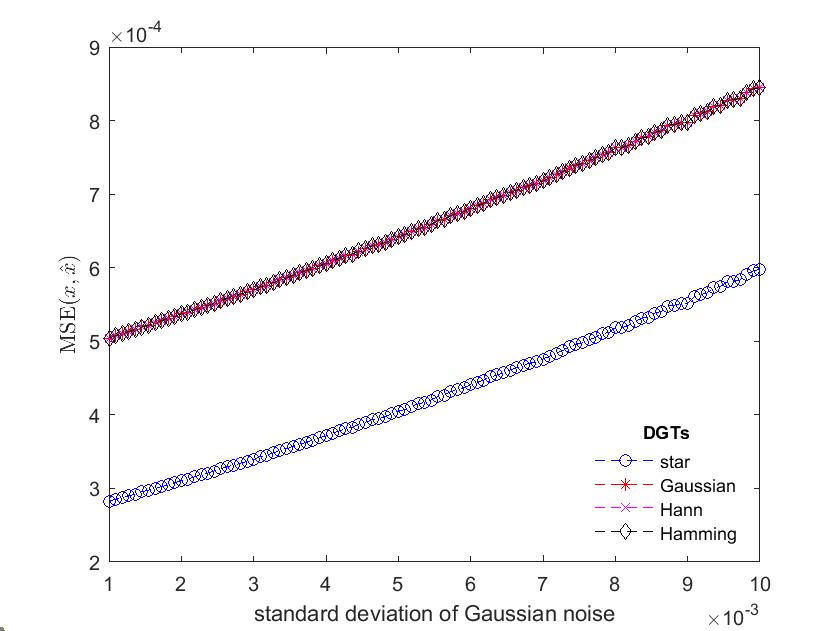}
  \end{minipage}
  }
  \caption{Rate of robustness for denoising real-world speech signals, with additive zero-mean Gaussian noise, having a standard deviation $\sigma$ taking values uniformly in $[0.001,0.01]$.}
\end{figure}

% \begin{table}[h]
%     \caption{Signals' details: top to bottom}
%     \centering
%     \scalebox{.8}{\begin{tabular}{|| c | c | c ||}
%     \hline
%      Labels & Samples & Ambient Dimension \& Lattice Parameters\\
%      \hline\hline
%     TwoChirp & 57 & $(L,a,b)=(57,1,19)$\\
%      \hline
%      Bumps & 45 & $(L,a,b)=(45,9,1)$\\
%      \hline
%      Cusp & 57 & $(L,a,b)=(57,1,19)$\\
%      \hline
%     \end{tabular}}
%     \label{sigdec}
% \end{table}

\begin{figure}
  \subfloat[]{%
  \begin{minipage}{\linewidth}
  \centering
  \includegraphics[width=0.8\textwidth]{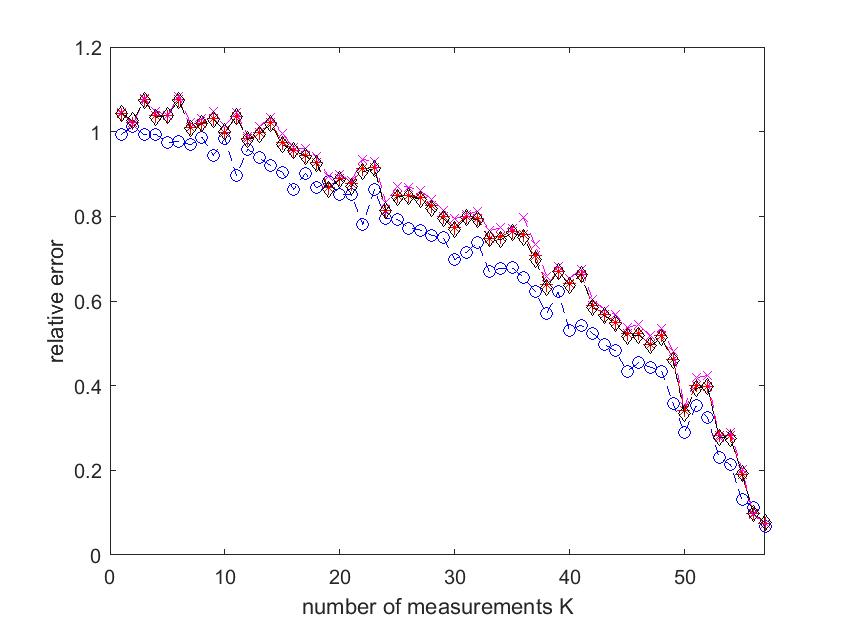}\hfill
  \includegraphics[width=0.8\textwidth]{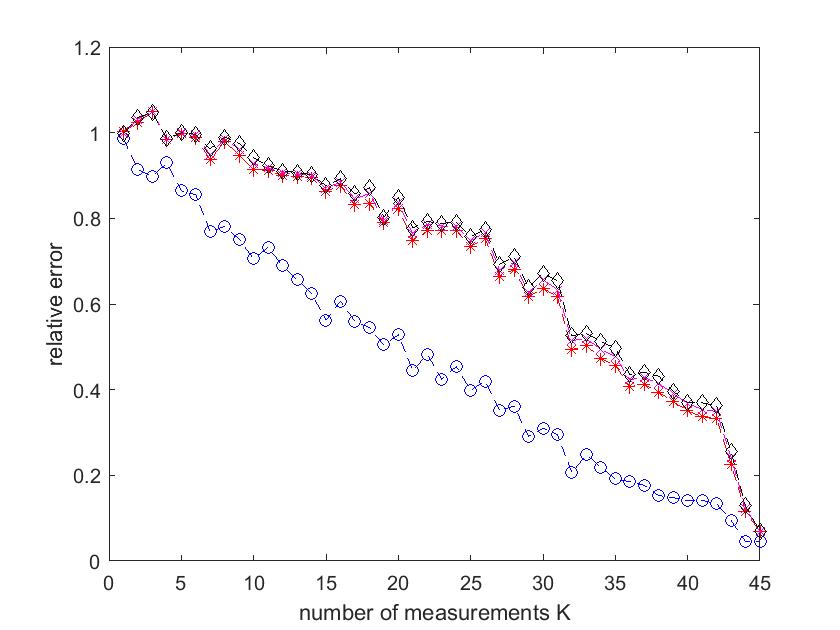}\hfill
  \includegraphics[width=0.8\textwidth]{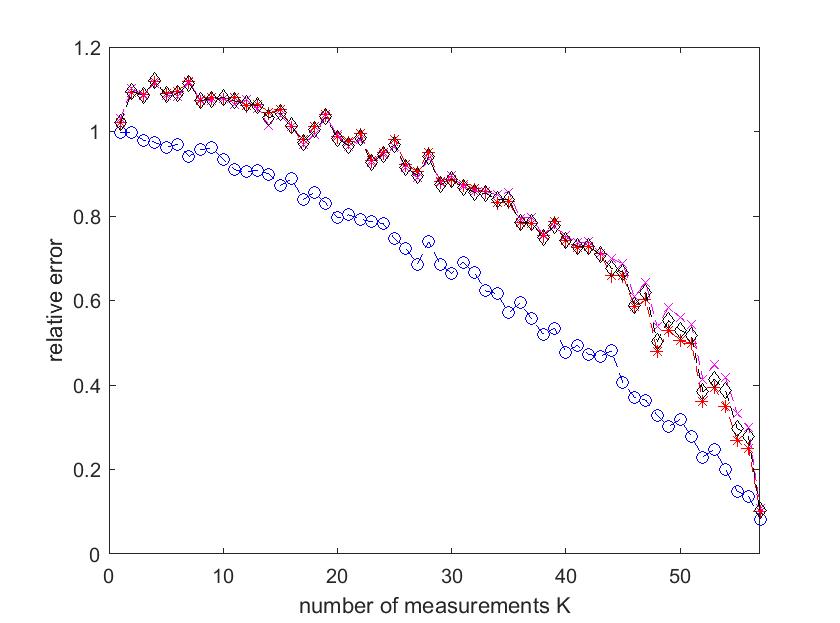}
  \end{minipage}
  }
  \caption{Rate of approximate success for CS with synthetic data, using a random Gaussian $A\in\mathbb{V}^{K\times L}$ and standard deviation of Gaussian noise added $\sigma=0.001$.}
  \label{}
\end{figure}

\begin{table}[h]
    \caption{Signals' details: top to bottom for both CS and denoising}
    \centering
    \scalebox{.75}{\begin{tabular}{|| c | c | c ||}
    \hline
     Labels & Samples & Ambient Dimension \& Lattice Parameters\\
     \hline\hline
     Denoising: 251-136532-0014 & 36240 & $(L,a,b)=(33915,51,19)$\\
     \hline
     Denoising: 3752-4944-0042 & 51360 & $(L,a,b)=(51051,21,23)$\\
     \hline
     Denoising: 5694-64038-0013 & 52880 & $(L,a,b)=(51051,33,17)$\\
     \hline
     CS: TwoChirp & 57 & $(L,a,b)=(57,1,19)$\\
     \hline
     CS: Bumps & 45 & $(L,a,b)=(45,9,1)$\\
     \hline
     CS: Cusp & 57 & $(L,a,b)=(57,1,19)$\\
     \hline
    \end{tabular}}
    \label{sigdec}
\end{table}

\bibliographystyle{IEEEtran}

\bibliography{mybib}

% that's all folks
\end{document}